\documentclass[a4paper,12pt]{article}

\usepackage{amsmath}
\usepackage{amsthm}
\usepackage{amssymb}  
\usepackage{latexsym} 
\usepackage[all]{xy}
\usepackage{comment}
\usepackage{rotating}
\usepackage{slashbox}
\usepackage{multirow}
\usepackage{rotating}
\usepackage{lscape}

\newcommand{\pp}{{\mathfrak {p}}}

\newcommand{\EE}{{\mathcal E}}
\newcommand{\CC}{{\mathcal C}}

\newcommand{\Mat}{\operatorname{Mat}}
\newcommand{\Char}{\operatorname{char}}

\newcommand{\Spec}{\operatorname{Spec}}

\newcommand{\diag}{\operatorname{diag}}
\newcommand{\id}{\operatorname{id}}
\newcommand{\mult}{\operatorname{mult}}

\newcommand{\GL}{\operatorname{GL}}

\newcommand{\GG}{{\mathcal G}_k}

\newcommand{\Kbar}{\overline{K}}
\newcommand{\kbar}{\overline{k}}

\newcommand{\summ}{\operatorname{sum}}
\newcommand{\isom}{ \cong }

\newcommand{\OK}{{\mathcal{O}_K}}

\newcommand{\PGL}{\operatorname{PGL}}

\newcommand{\Gal}{\operatorname{Gal}}

\newcommand{\Aff}{{\mathbb A}}
\newcommand{\PP}{{\mathbb P}}
\newcommand{\Q}{{\mathbb Q}}

\newcommand{\Norm}{\operatorname{Norm}}

\newcommand{\SL}{\operatorname{SL}}

\newcommand{\Proj}{\operatorname{Proj}}

\newcommand{\Z}{{\mathbb Z}}
\newcommand{\Aut}{\operatorname{Aut}}

\newfont{\wncyr}{wncyr10 at 12pt}
\newfont{\wncyrten}{wncyr10 at 10pt}

\newenvironment{Proof}{\par\noindent{\sc Proof:}}%
                      {\hspace*{\fill}\nobreak$\Box$\par\medskip}
\newenvironment{ProofOf}[1]{\par\noindent{\sc Proof of #1:}}%
                       {\hspace*{\fill}\nobreak$\Box$\par\medskip}
\newenvironment{myitemize}
{\begin{itemize}
\setlength{\itemsep}{1pt}
\setlength{\parskip}{0pt}
\setlength{\parsep}{0pt}}
{\end{itemize}}

\newtheorem{Proposition}{Proposition}[section]
\newtheorem{Theorem}[Proposition]{Theorem}
\newtheorem{Lemma}[Proposition]{Lemma}
\newtheorem{Corollary}[Proposition]{Corollary}

\theoremstyle{definition}
\newtheorem{Definition}[Proposition]{Definition}

\newtheorem{Example}[Proposition]{Example}

\newcounter{nootje}
\begin{document}

\title{Counting Models of Genus One Curves}
\author{Mohammad Sadek}
\date{}
\maketitle
\begin{abstract}{\footnotesize
Let $C$ be a soluble smooth genus one curve over a Henselian discrete valuation field. There is a unique minimal Weierstrass equation defining $C$ up to isomorphism. In this paper we consider genus one equations of degree $n$ defining $C$, namely a (generalised) binary quartic when $n=2$, a ternary cubic when $n=3$, and a pair of quaternary quadrics when $n=4$. In general, minimal genus one equations of degree $n$ are not unique up to isomorphism. We explain how the number of these equations varies according to the Kodaira symbol of the Jacobian of $C$. Then we count these equations up to isomorphism over a number field of class number $1$.}
\end{abstract}

\section{Introduction}
\label{sec:intro}

Let $E$ be an elliptic curve defined over a number field $K$. It is known that $E(K)$ is a finitely generated abelian group. If we manage to find $E(K)/nE(K)$ for any integer $n\ge 2$, then a finite amount of computation will determine $E(K)$. The method of $n$-descent is one of the methods which enable us to get a bound on $E(K)/nE(K),\;n\ge2.$ Indeed, the $n$-descent computes the $n$-Selmer group of $E$ which contains $E(K)/nE(K)$ as a subgroup. The difference between the two groups is the $n$-torsion of the Tate-Shafarevich group.

 An element of the $n$-Selmer group can be represented as a geometric object, namely as an everywhere locally soluble curve defined by a {\em genus one equation of degree $n$} given as follows:

 If $n=1,$ a Weierstrass equation \begin{equation} y^2 + a_1 x y + a_3 y = x^3 + a_2 x^2 + a_4 x + a_6.\end{equation}
Two genus one equations of degree $1$ with coefficients in a Dedekind domain $R$ are {\em $R$-equivalent} if they are related by the substitutions:
$x' = u^2x + r$ and $y' = u^3y + su^2x + t,$ where $r,s,t\in R,\;u\in R^*.$ The group of transformations $[u;r,s,t]$ is called $\mathcal{G}_1(R)$. We set $\det([u;r,s,t])=u^{-1}.$

 If $n=2,$ a (generalised) binary quartic \begin{equation}y^2 + (\alpha_0 x^2 + \alpha_1 xz + \alpha_2z^2) y = a x^4 + b x^3z + c x^2z^2 + d xz^3 + ez^4.\end{equation}Two genus one equations of degree $2$ with coefficients in $R$ are {\em $R$-equivalent} if they are related by the substitutions: $x'=m_{11}x+m_{21}z,\;z'=m_{12}x+m_{22}z$ and $y'=\mu^{-1}y+r_0x^2+r_1xz+r_2z^2,$ where $\mu\in R^*,\; r_i\in R$ and $M=(m_{ij})\in \GL_2(R)$. The group of transformations $[\mu,(r_i),M]$ is called $\mathcal{G}_2(R)$. We set $\det([\mu,(r_i),M])=\mu\det M.$

 If $n=3,$ a ternary cubic \begin{equation}F(x,y,z)=a x^3 + b y^3 + c z^3 + a_2x^2 y + a_3x^2 z + b_1y^2x+ b_3 y^2z + c_1z^2x + c_2 z^2y + m x y z=0.\end{equation} Two genus one equations of degree $3$ with coefficients in $R$ are {\em $R$-equivalent} if they are related by multiplying by $\mu\in R^*$, then substituting
 $x'=m_{11}x+m_{21}y+m_{31}z,\;y'=m_{12}x+m_{22}y+m_{32}z$ and $z'=m_{13}x+m_{23}y+m_{33}z$, where $M=(m_{ij})\in \GL_3(R).$
 The group of transformations $[\mu, M]$ is called $\mathcal{G}_3(R)$. Set $\det([\mu,M])=\mu\det M.$

 If $n=4$, then $C$ is described by $2$ quadratic forms in $4$ variables. Two genus one equations of degree $4$ with coefficients in $R$ are {\em $R$-equivalent} if they are related by the substitutions: $F'_i=m_{i1}F_1+m_{i2}F_2$ for $M=(m_{ij})\in\GL_2(R),$
and then $x'_j=\sum_{i=1}^4n_{ij}x_i$ for $N=(n_{ij})\in\GL_4(R).$
The group of transformations $[M,N]$ is called $\mathcal{G}_4(R)$. Set $\det([M,N])=\det M\det N.$

For $n\le4$, we associate invariants $c_{4,\phi}$, $c_{6,\phi}$ and {\em discriminant} $\Delta_{\phi}$ to a genus one equation $\phi$ of degree $n$
where $\Delta_{\phi}=(c_{4,\phi}^3-c_{6,\phi}^2)/1728.$ Moreover, $\phi$ defines a smooth curve of genus one if and only if $\Delta_{\phi}\ne 0.$ The invariants $c_{4,\phi},c_{6,\phi}$ and $\Delta_{\phi}$ are of weights $r=4,6$ and $12$ respectively, in other words if $F\in\{c_{4,\phi},c_{6,\phi},\Delta_{\phi}\}$, then $F\circ g= (\det g)^rF$ for all $g\in\mathcal{G}_n$.
For more details on these genus one equations and their classical invariants see \cite{AKMMMP}. We scale the above invariants according to \cite{FiInvariants} and \cite{FiStCr}.

For using $n$-descent to search for points on an elliptic curve $E,$ we need the coefficients of the genus one equations described above to be small. That can be achieved by {\em Reduction} and {\em Minimisation}. By reducing genus one equations, we mean reducing the size of the coefficients by a unimodular linear change of coordinates, which does not change the invariants. To minimise genus one equations, we need to make the associated invariants smaller. For references see \cite{swinnerton} and \cite{CremonaStollMini} for $n=2,$ \cite{Fiminimise} for $n=3,$ and \cite{WomackThesis} for $n=4.$ An algorithmic approach to the minimisation problem can be found in \cite{FiStCr}. In this paper we will be concerned with the minimisation question.

We work throughout over a perfect Henselian discrete valuation field $K$ with algebraic closure $\Kbar$ and ring of integers $\OK$. We fix a uniformiser $t$, a normalisation $\nu$ and write $k=\OK/t\OK$. Set $S=\Spec\OK$.

 \begin{Definition}\label{def:integral-minimal}
 A genus one equation $\phi$ of degree $n$, $n\le4$, with $\Delta_{\phi}\ne 0$ is
 \begin{myitemize}
\item[(a)] {\em integral} if the defining polynomials have coefficients in $\OK.$
\item[(b)] {\em minimal} if it is integral and $\nu(\Delta_{\phi})$ is
minimal among all the valuations of the discriminants of the integral genus one equations $K$-equivalent to $\phi.$
\end{myitemize}
\end{Definition}

It has been shown in \cite{FiStCr} and \cite{SadekMinimal} that if $C$ is a soluble smooth curve defined by a genus one equation $\phi$ of degree $n$, then the minimal discriminant associated to $\phi$ has the same valuation as the minimal discriminant associated to the Jacobian elliptic curve.

The $S$-scheme $\CC$ defined by an integral genus one equation $\phi:y^2+g(x,z)y=f(x,z)$ of degree $2$
is the scheme obtained by glueing
$ \{ y^2+g(x,1)y = f(x,1) \} \subset \Aff_S^2 $ and
$ \{ v^2+g(1,u)v = f(1,u) \} \subset \Aff_S^2 $ via $x=1/u$ and $y=x^2v$. It comes with
a natural morphism $\CC \to \PP^1_S$ given on these affine pieces
by $(x,y) \mapsto (x:1)$ and $(u,v) \mapsto (1:u)$.
The $S$-scheme defined by an integral genus one equation $\phi$ of degree $n$, $n=1$ or $n\ge3$, is
the subscheme $\CC \subset \PP^{m}_S$ defined by $\phi$, where $m=2$ when $n=1$, and $m=n-1$ when $n\ge3$. In \cite{SadekMinimal}, we gave the following criterion for minimality.
  \begin{Theorem}
\label{thm:minimality}
Let $\phi$ be an integral genus one equation of degree $n,\;n\le4$. Assume moreover that $\phi$ defines a normal $S$-scheme $\CC$, and that $\CC_K$ is smooth with $\CC_K(K)\ne\emptyset.$
Then $\phi$ is minimal if and only if the minimal desingularisation of $\CC$ is isomorphic to the minimal proper regular model of the Jacobian of $\CC_K$.
\end{Theorem}

 Unlike elliptic curves, the minimisations of integral genus one equations of degree $n$ are not unique under the action of the group $\mathcal{G}_n(\OK)$. In this paper we study the non-uniqueness of minimisation. In order to achieve this goal we generalise the scheme-theoretic  terminology of minimal Weierstrass models for elliptic curves to {\em minimal degree-$n$-models} for genus one curves.

\begin{Definition}
\label{def:Weierstrass model}
Let $\phi$ be a genus one equation of degree $n$. Let $C\to \PP_K^{m}$ be the genus one curve over $K$ defined by $\phi,$ where $m=2$ when $n=1$ and $m=n-1$ otherwise.
A {\em (minimal) degree-$n$-model} for $C\to\PP^m_K$ is a pair $(\CC,\alpha)$ where $\CC\subset \PP^m_S$ is an
$S$-scheme defined by a (minimal) integral genus one equation of degree $n$, and
$\alpha:\CC_K\isom C$ is an isomorphism defined by a $K$-equivalence of genus one equations of degree $n,$ i.e., $\alpha$ is defined by an element of $\mathcal{G}_n(K)$. If there is no confusion, then we will omit mentioning $\alpha$ in $(\CC,\alpha)$ and write $\CC$ instead.

 {\em An isomorphism} $\beta:(\CC_1,\alpha_1)\isom(\CC_2,\alpha_2)$ of degree-$n$-models is an isomorphism $\beta:\CC_1\isom\CC_2$ of $S$-schemes defined by an $\OK$-equivalence of genus one equations of degree $n,$ i.e., $\beta$ is defined by an element of $\mathcal{G}_n(\OK),$ with $\beta_K = \alpha_2^{-1} \alpha_1$.
  \end{Definition}

  Now we show why we insist on remembering the isomorphism $\alpha$ in our definition of a degree-$n$-model. In the case $C(K)\ne\emptyset$, we pick $P\in C(K)$. Then we identify $(C,P)$ with the Jacobian elliptic curve $E$. The automorphism group $\Aut(C)$ of $C$ fits in an exact sequence
 $$0\to E\to\Aut(C)\to\Aut(E,0)\to 0.$$
 The first map is $Q\mapsto\tau_Q$, where $\tau_Q$ is the translation by $Q$. Let $H$ be a hyperplane section on $C$. We are interested in automorphisms $\lambda$ of $C$ such that $\lambda^*H\sim H$. But $\tau_Q^*H\sim H$ if and only if $nQ=0$. Hence the elements of $\PGL_n(K)$ that act fixed-point-free on $C$ correspond precisely to $E[n](K)$.
 Therefore, we can have more than one $K$-isomorphism between $\CC_K$ and $C$ defined by elements in $\mathcal{G}_n(K)$.

 Liu proved that there is a bijection between minimal degree-$2$-models for $C\to\PP^1_K$, up to isomorphism, and the multiplicity-$1$ irreducible components in the special fiber of the quotient of the minimal proper regular model by a hyperelliptic involution, see (\cite{LiuModeles}, Corollaire 5).

 In \S \ref{sec:computing in E min}, we perform explicit computations on the minimal proper regular model of an elliptic curve. In \S \ref{sec:counting minimal degree n models}, these computations are used to count minimal degree-$n$-models for a soluble smooth genus one curve $C\to\PP^{n-1}_K$, $n\in\{2,3,4\}$, up to isomorphism. We relate the number of minimal degree-$n$-models for $C\to\PP_{K}^{n-1}$, up to isomorphism, to the cardinality of a finite set which depends only on the Kodaira Symbol of the Jacobian of $C$ and on the divisor defining the morphism $C\to\PP^{n-1}_K$, see Theorem \ref{thm:counting models}.
 Then, using strong approximation, we prove that the number of minimal global degree-$n$-models for $C\to\PP^{n-1}_K$ over a number field $K$ of class number 1 is the product of minimal degree-$n$-models for $C\to\PP^{n-1}_{K_{\nu}}$, where $K_{\nu}$ is the completion of $K$ at a finite place $\nu$ and $\nu$ runs over a finite set of bad primes of the Jacobian of $C$.

We have to mention that over a Henselian discrete valuation field, our counting results work for any residue characteristic if the Jacobian has multiplicative reduction, while for additive reduction we assume $\Char k\ne 2$ when $n=2$, and $\Char k\ne 2,3$ when $n=3,4$. It seems to us that these conditions on the residue characteristic for the additive reduction case are not necessary, but we have not been able to eliminate them.

\section{Isomorphisms of degree-$n$-models}
\label{sec:isomorphisms of degree n models}

Let $C$ be a smooth genus one curve over $K$. Assume that $C(K)\ne \emptyset$. Let $E$ be the Jacobian of $C$ with minimal proper regular model $E^{min}.$ We define an {\em $S$-model} for $C$ to be an integral, projective, flat, normal $S$-scheme $f:X\to S$ of dimension 2 such that its generic fiber $\CC_K$ is isomorphic to $C$. The special fiber of $\CC$ will be called $\CC_k$.

Let $P\in C(K).$ Fix an isomorphism $\beta:C\to E$ such that $\beta(P)=0_E.$ The isomorphism $\beta$ identifies $(C,P)$ with $(E,0_E).$ We observe that $\beta$ extends to an isomorphism between the minimal proper regular model of $C$ and $E^{min}.$
 In \cite{SadekMinimal} we proved that a minimal degree-$n$-model $\CC$ for a smooth genus one curve $C\to\PP^{n-1}_K$ is normal, and hence is an $S$-model for $C$. Moreover, we proved that the minimal desingularisation of $\CC$ is isomorphic to $E^{min}$, see Theorem \ref{thm:minimality}. Furthermore, we have
$$\CC=\Proj(\bigoplus_{m=0}^{\infty}H^0(E^{min},\mathcal{O}_{E^{min}}(mD))),$$
where $D$ is a divisor on $E^{min}$ which intersects the irreducible components of the strict transform of $\CC_k$ in $E^{min}_k$ and intersects no other components in $E^{min}_k$. The divisor $D|_E$ is linearly equivalent to the hyperplane section divisor defining $C$.

Let $\Gamma$ be an irreducible component of $\CC_k.$ By the {\em type} of $\Gamma$ we mean the ordered pair $(\mult_k(\Gamma),\deg_k(\Gamma))$, where $\mult_k(\Gamma)$ and $\deg_k(\Gamma)$ are the multiplicity and degree of $\Gamma$ respectively.

\begin{Theorem}
\label{thm:injectivity2}
Let $C$ be a smooth genus one curve with $C(K)\ne\emptyset$. Let $(\CC_1,\alpha_1)$ and $(\CC_2,\alpha_2)$ be two minimal degree-$n$-models for $C\to\PP^{n-1}_K$, $n=2,3,4$. Set $\alpha=\alpha_2^{-1}\alpha_1$ and denote its extension $\CC_1\dashrightarrow\CC_2$ by $\tilde{\alpha}.$ Then the following statements are equivalent.
 \begin{myitemize}
 \item[(i)] $(\CC_1,\alpha_1)$ and $(\CC_2,\alpha_2)$ are isomorphic as degree-$n$-models for $C\to\PP^{n-1}_K.$
  \item[(ii)] For every irreducible component $\Gamma$ of $(\CC_2)_k,$ $\tilde{\alpha}^*\Gamma$ is an irreducible component of $(\CC_1)_k$ with the same type as $\Gamma$.
  \end{myitemize}
\end{Theorem}
\begin{Proof}
$(i)$ implies $(ii)$: Let $\alpha:=\alpha_2^{-1}\alpha_1:(\CC_1)_K\to(\CC_2)_K.$ By the definition of isomorphic degree-$n$-models, the map $\alpha$ extends to an $S$-isomorphism $\tilde{\alpha}:\CC_1\to\CC_2$ which is defined by an element in $\mathcal{G}_n(\OK).$ Therefore, if $\Gamma$ is an irreducible component of $(\CC_2)_k,$ then $\tilde{\alpha}^*\Gamma$ is an irreducible component of $(\CC_1)_k$ with the same type as $\Gamma.$

$(ii)\Rightarrow(i)$: We want to show that $\tilde{\alpha}$ is defined by an element in $\mathcal{G}_n(\OK).$
 Since the minimal proper regular models of $(\CC_1)_K$ and $(\CC_2)_K$ are identified with $E^{min}$, statement $(ii)$ implies that both $(\CC_1)_k$ and $(\CC_2)_k$ have the same strict transform in $E^{min}_k.$
Let $D_i$ be a defining divisor of $\CC_i$ as a contraction in $E^{min},$ i.e.,
$$\CC_i=\Proj(\bigoplus_{m\ge0}H^0(E^{min},\mathcal{O}_{E^{min}}(m D_i))).$$
Set $\mathcal{L}_i=\mathcal{O}_{E^{min}}(D_i)$ and $\mathcal{L}_{i,k}=\mathcal{L}_i|_{E^{min}_k}.$
Let $\Gamma$ be an irreducible component of $(\CC_2)_k.$ Consider the strict transform $\tilde{\Gamma}$ of the irreducible components $\Gamma$ and $\tilde{\alpha}^*\Gamma$ in $E^{min}_k.$ Since $\Gamma$ and $\tilde{\alpha}^*\Gamma$ have the same type, it follows that
  $\deg_k \mathcal{L}_{1,k}|_{\tilde{\Gamma}}=\deg_k \mathcal{L}_{2,k}|_{\tilde{\Gamma}}.$ For any irreducible component $\Lambda$ which is not a strict transform of a component of $(\CC_i)_k,$ we have $\deg_k \mathcal{L}_{i,k}|_{\Lambda}=0.$ Now since each irreducible component of $E^{min}_k$ is isomorphic to $\PP^1_k$, we have $D_1|_{E^{min}_k}$ is linearly equivalent to $D_2|_{E^{min}_k}$ as they have the same degree on every irreducible component of $E^{min}_k$. It follows that $\mathcal{L}_{1,k}\isom \mathcal{L}_{2,k}.$

   We are given that $C(K)\ne\emptyset$, $\mathcal{L}_1|_{(\CC_1)_K}\isom\mathcal{L}_2|_{(\CC_2)_K}$ and we established the isomorphism $\mathcal{L}_{1,k}\isom\mathcal{L}_{2,k}.$ These conditions imply that $\mathcal{L}_1\isom\mathcal{L}_2,$ see (\cite{Liubook}, Exercise 9.1.13 (b)). Therefore, $H^0(E^{min},\mathcal{L}_1)$ and $H^0(E^{min},\mathcal{L}_2)$ are isomorphic as $\OK$-modules, and $\tilde{\alpha}$ is a change of basis of a free $\OK$-module of rank $n$. Hence $\tilde{\alpha}$ is defined by an element in $\mathcal{G}_n(\OK).$
\end{Proof}

\begin{Corollary}
\label{cor:unique model when IO or I1}
Let $C$ be a smooth genus one curve. Assume that $C(K)\ne\emptyset$. Assume moreover that the Jacobian $E$ of $C$ has either reduction types ${\rm I_0}$ or ${\rm I_1}$. Then there is a unique minimal degree-$n$-model for $C\to\PP^{n-1}_K$.
\end{Corollary}

\section{Computing in $E^{min}$}
\label{sec:computing in E min}

The following lemma will play an essential rule in the way we construct divisors on minimal proper regular models.

 \begin{Proposition}
\label{prop:auxiliary divisors}
 Let $C/K$ be a smooth genus one curve with $C(K)\ne\emptyset$. Let $E$ be the Jacobian of $C$ with minimal proper regular model $E^{min}$.  Fix a closed point $x\in E^{min}_k$ such that $x$ lies on one and only one irreducible component $\Gamma$ of $E^{min}_k$, of multiplicity $r\ge 1$. Then there exists a point $P\in C$ such that $\overline{\{P\}}\cap E^{min}_k=\{x\}$ and $[K(P):K]=r,$ where $\overline{\{P\}}$ is the Zariski closure of $\{P\}$ in $E^{min}$.
\end{Proposition}
\begin{Proof}
See (\cite{Liubook}, Exercise 9.2.11 (c)).
\end{Proof}

It is easy to check that
if $L$ is a totally ramified extension of $K$ with $[L:K]=m,\;m=2,3,4,$ then the maximal unramified extension $L^{un}$ of $L$ is a Galois extension of the maximal unramified extension $K^{un}$ of $K$ with $L^{un}=K^{un}(\sqrt[m]{t})$.

Let $E$ be an elliptic curve over $K,$ with minimal proper regular model $E^{min}.$ Let $E^0(K)$ be the group of rational points of $E$ with non-singular reduction. In this section we define a map $\delta_m:\Phi_K^m(E)\to \Phi_K(E),$ where $\Phi_K^m(E)$ is the set of multiplicity-$m$ irreducible components of $E^{min}_k,$ and $\Phi_K(E):=\Phi_K^1(E)$ is the group of components $E(K)/E^0(K).$ In the remainder of this section we assume that the residue field $k$ is algebraically closed and hence every finite extension of $K$ is a totally ramified Galois extension, and the reduction type of $E$ is split.

\subsection{The component group}
\label{sec:the components group}

 The component group $\Phi_K(E)=E(K)/E^0(K)$ is finite.
  More precisely, if $E$ has multiplicative reduction, then $E(K)/E^0(K)$ is a cyclic group of order $\nu(\Delta);$ otherwise, $E(K)/E^0(K)$ has order $1,2,3,$ or $4$, see (\cite{Sil2}, Chapter IV, Corollary 9.2).

 If $E$ has reduction $I_n,\;n\ge0,$ then $E^{min}_{k}$ consists of an $n$-gon. There exists an isomorphism $\delta_1:\Phi_{K}(E)\isom\mathbb{Z}/n\mathbb{Z}$ such that if we number the irreducible components in $\Phi_K(E)$ consecutively from $0$ to $n-1$, then the image of the $i$th component in $\Phi_{K}(E)$ is $i\in\Z/n\Z$.
 If $E$ has additive reduction, then we fix an isomorphism $\delta_1:\Phi_{K}(E)\isom \mathbb{Z}/2\mathbb{Z}\times\mathbb{Z}/2\mathbb{Z}$ when $E$ has reduction of type ${\rm I}_{2m}^*,\;m\ge0,$ and $\delta_1:\Phi_K(E)\isom\mathbb{Z}/n\mathbb{Z}$ for some $n\in\{1,2,3,4\}$ for the other additive reduction types. From now on we will not distinguish between $\Phi_K(E)$ and its isomorphic image under $\delta_1$.

We note that if $E$ has one of the reduction types ${\rm I}_n,n\ge0,{\rm II},{\rm III},$ or ${\rm IV},$ then $E^{min}_k$ consists only of multiplicity-$1$ components, hence $\Phi_K^m(E)=\emptyset$ when $m\ge2$.

  We define the map $\delta_m:\Phi_K^m(E)\to\Phi_K(E),\;m\ge2,$ as follows: Every irreducible component in $\Phi_K^m(E)$ is isomorphic to $\PP_k^1$. Let $x$ be a $k$-point of $\Theta\in\Phi_K^m(E)$ which lies on no other component. Let $P\in E(\Kbar)$ be a point which reduces to $x$ such that $[K(P):K]=m,$ see Proposition \ref{prop:auxiliary divisors}. Let $\sigma$ be a generator of $\Gal(K(P)/K)$. The sum of the Galois orbit $P+\ldots+P^{\sigma^{m-1}}$ of $P$ is a point in $E(K).$ We set $\delta_m(\Theta)$ to be the image of this sum in $\Phi_K(E)=E(K)/E^0(K).$ It will be clear from the description of the map $\delta_m$ given below that $\delta_m(\Theta)$ does not depend on the choice of $x$ nor on the choice of $P.$

Since $K(P)/K$ is a totally ramified extension, the Galois group $\Gal(K(P)/K)$ is the inertia group of $K(P)/K.$ In other words, $\Gal(K(P)/K)$ fixes every irreducible component in $E^{min}_k.$ Hence if $P\in E(\Kbar)$ lies above an irreducible component $\Gamma$ of $E^{min}_k$, then every point in the Galois orbit of $P$ lies above $\Gamma$ as well.

 \subsection{Reduction type ${\rm I}_n^*,\;n\ge0$}
 \label{subsec:In*}

Assume that $E/K$ has reduction of type ${\rm I}_n^*,\;n\ge0.$ The special fiber $E^{min}_k$ contains a sequence of $(n+1)$ multiplicity-$2$ components and no components of higher multiplicities. Tate's algorithm, see (\cite{Sil2}, Chapter IV, \S9), can be used to write down explicit equations for the components in $\Phi_K^2(E)$, hence we can  determine conditions for a point in $E$ defined over $K(\sqrt{t})$ to lie above one of these multiplicity-$2$ components, and compute the sum of the Galois orbit of this point.
We will write $a_{i,r}$ for $t^{-r}a_i,$ and $x_r,y_r$ for $t^{-r}x,t^{-r}y$ respectively.

 The proof of Tate's algorithm shows that if $E$ has additive reduction of one of the types ${\rm I}_n^*,\;n\ge0,\;{\rm IV}^*,\;{\rm III}^*$ or ${\rm II}^*$, then we can assume that there exists a minimal Weierstrass equation of the form
 $$y^2z+a_1xyz+a_3yz^2=x^3+a_2x^2z+a_4xz^2+a_6z^3,\;a_i\in\OK,$$
 for $E/K,$ such that $t\mid a_1,a_2,t^2\mid a_3,a_4$ and $t^3\mid a_6.$ Moreover, the identity component is given by $\Gamma_0:z=0$ and is attached to a multiplicity-$2$ component $V_{-1}:y_1^2=0$.

  If $E$ has reduction of type ${\rm I}_n^*,\;n\ge0,$ then the non-identity multiplicity-$1$ component attached to $V_{-1}$ is $\Gamma:x_1+\tilde{a}_{2,1}=0$. The other multiplicity-$2$ components of $E^{min}$ are given by
\begin{displaymath}
V_l: \left\{ \begin{array}{ll}
y_u^2+\tilde{a}_{3,u}y_u-\tilde{a}_{6,2u}=0 & \textrm{ if } l=2u-3\textrm{ is odd,}\\
\tilde{a}_{2,1}x_u^2+\tilde{a}_{4,u+1}x_u+\tilde{a}_{6,2u+1}=0& \textrm{ if }l=2u-2\textrm{ is even.}
\end{array} \right.
\end{displaymath}
In addition, $V_l$ consists of two distinct lines precisely when $l=n=\nu(\Delta)-6.$
The special fiber $E^{min}_k$ is as in the following figure. \\

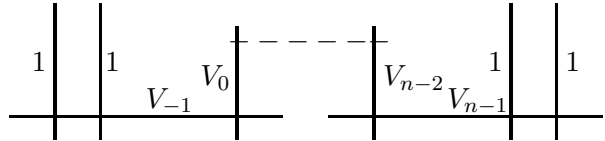
\setlength{\unitlength}{0.6cm}\linethickness{.3mm}
\small
\begin{picture}(18,5)
\put(5,2){\line(1,0){6}}\put(8,2.2){$V_{-1}$}
\put(6,1.5){\line(0,1){3}}\put(5.5,3){$1$}
\put(7,1.5){\line(0,1){3}}\put(7.1,3){$1$}
\put(10,1.5){\line(0,1){2.5}}\put(9.2,2.7){$V_0$}
\put(9.8,3.5){$------$}
\put(13,1.5){\line(0,1){2.5}}\put(13.2,2.7){$V_{n-2}$}
\put(17,1.5){\line(0,1){3}}\put(17.2,3){1}
\put(16,1.5){\line(0,1){3}}\put(15.5,3){1}
\put(12,2){\line(1,0){6}}\put(14.6,2.2){$V_{n-1}$}
\end{picture}\\
Let $K_2=K(\sqrt{t})$ be the unique quadratic tame Galois extension of $K.$ A point lying above a $k$-point on $V_l,$ where $l=2u-2$, and on no other component is a point in $ E(K_2)$ of the form
 $(\alpha t^{u+1}+\beta t^{u+1/2},y),$ where $y\in K_2,\;t^{u+1}\mid y$, $\alpha\in\OK$ and $\beta\in \mathcal{O}_K^*$. We have $\beta\in \mathcal{O}^*_K,$ since otherwise this point would reduce to a point on $V_{l+1}.$ Similarly, a point lying above a $k$-point on $V_l,$ where $l=2u-3$, and on no other component is a point in $ E(K_2)$ of the form
 $(x,\alpha t^{u+1}+\beta t^{u+1/2}),$ where $x\in K_2,\;t^u\mid x$, $\alpha\in\OK$ and $\beta\in \mathcal{O}_K^*$.

 We recall the addition formula on elliptic curves. Let $P_1,P_2\in E(\Kbar).$ If $P_1\ne\pm P_2,$ then
 $$x(P_1+P_2)=\lambda^2+a_1\lambda-a_2-x(P_1)-x(P_2),\textrm{ and, }
y(P_1+P_2)=-(\lambda+a_1) x(P_1+P_2)-\mu-a_3,$$
where $$\lambda=\frac{y(P_2)-y(P_1)}{x(P_2)-x(P_1)},\;\mu=\frac{y(P_1)x(P_2)-y(P_2)x(P_1)}{x(P_2)-x(P_1)}.$$

   \begin{Lemma}
   \label{lem:y2=0, x3=0}
   Assume $E$ has one of the reduction types ${\rm I}_n^*,\;n\ge0,\;{\rm IV}^*,\;{\rm III}^*,$ or ${\rm II}^*$. Let $K_2=K(\sqrt{t})$ with $\Gal(K_2/K)=\langle\sigma\rangle$. Then a point $P\in E(K_2)$ lying above the component $y_1^2=0$ satisfies $P+P^{\sigma}\in E^0(K)$.
   \end{Lemma}
   \begin{Proof}
   If $P^{\sigma}=-P,$ then $P+P^{\sigma}=0\in E^0(K)$ and we are done. So, assume $x(P)\not\in K.$ Now $P=(\alpha't+\beta't^{3/2},\alpha t^{2}+\beta t^{3/2})$, where $\alpha',\beta'\in \OK$, $\alpha\in\OK$ and $\beta\in \mathcal{O}_K^*.$ In the addition formula given above, $\lambda=\beta/\beta'\in K.$
  Since $\nu(\beta')\ge 0,$ it follows that $\nu(x(P+P^{\sigma}))\le0$ because $t\mid a_1,a_2$. Hence the reduction of $P+P^{\sigma}$ is not a singular point, i.e., $P+P^{\sigma}\in E^0(K).$
   \end{Proof}

 \begin{Proposition}
 \label{prop:sum of points In*}
 Assume that $E/K$ has reduction type ${\rm I}_n^*,\;n\ge0.$ Let $K_2=K(\sqrt{t})$ with $\Gal(K_2/K)=\langle\sigma\rangle.$
 \begin{myitemize}
 \item[(i)] The multiplicity-$1$ component $\Gamma:x_1+\tilde{a}_{2,1}=0$ is of order $2$ in $\Phi_K(E)$.
   \item[(ii)] Let $P\in E(K_2)$ be a point lying above $V_l$ where $l=2u-3$. Then $P+P^\sigma\in E^0(K).$
  \item[(iii)] Let $P\in E(K_2)$ be a point lying above $V_l$ where $l=2u-2$. Then $P+P^\sigma$ reduces to a point on $\Gamma.$
\end{myitemize}
 \end{Proposition}
 \begin{Proof}
 $(i)$  If $P$ lies above $\Gamma:x_1+\tilde{a}_{2,1}=0,$ then $-P$ lies above $\Gamma$ as well. It follows that the inverse of $\Gamma$ when considered as an element of $\Phi_K(E)$ is itself.

 $(ii)$ If $P^{\sigma}=-P,$ then $P+P^{\sigma}=0\in E^0(K)$ and we are done. So, assume that $x(P)\not\in K$. We can assume that $P=(\alpha't^{u}+\beta't^{u+1/2},\alpha t^{u+1}+\beta t^{u+1/2})$, where $\alpha,\beta\in \OK$ with $\nu(\beta)=0$ and $\alpha',\beta'\in \OK$. We have
 $\lambda=\beta/\beta'$, hence $\nu(\lambda)\le 0$. Thus $\nu(x(P+P^{\sigma}))\le 0,$ and $P+P^{\sigma}\in E^0(K).$

 $(iii)$ We can assume that $P=(\alpha t^{u+1}+\beta t^{u+1/2},y)$, where $y\in K_2,\;t^{u+1}\mid y$ and $\alpha,\beta\in \OK$ with $\nu(\beta)= 0.$ Therefore, we have $\nu(\lambda)>0$ where $\lambda:=(y(P)-y(P^{\sigma}))/(x(P)-x(P^{\sigma})).$ Following the addition formula we have $x(P+P^{\sigma})=\lambda^2+a_1\lambda-a_2-x(P)-x(P^{\sigma}).$ Dividing by $t$ and using that $t\mid a_1,a_2$, we get that $x_1(P+P^{\sigma})=-\tilde{a}_{2,1}$ mod $t$.
\end{Proof}

  In Proposition \ref{prop:sum of points In*}, since $\Gamma:x_1+\tilde{a}_{2,1}=0$ has order 2, we have $\delta_1(\Gamma)=2\in\Phi_K(E)\isom \mathbb{Z}/4\mathbb{Z},$ when $n$ is odd. When $n$ is even, we know that $\delta_1:\Phi_K(E)\isom \mathbb{Z}/2\mathbb{Z}\times \mathbb{Z}/2\mathbb{Z},$ therefore every non-identity irreducible component has order $2$. From now on, we will fix $\delta_1$ such that $\delta_1(\Gamma)=(1,1).$

 \begin{Corollary}
 \label{cor:sum In*}
  Assume $E/K$ has reduction type ${\rm I}_n^*,\;n\ge0.$ Let
  $\Phi_K^2(E)=\{V_{-1},V_0,\ldots,V_{n-1}\}$
   be as above. Then $\delta_2:\Phi_K^2(E)\to\Phi_K(E)$ is determined according to the following table.

   \begin{center}
 \begin{tabular}{|c|c|c|c|c|}
  \hline
 $(n,l)$&$(2m+1,2u-3)$&$(2m+1,2u-2)$&$(2m,2u-3)$&$(2m,2u-2)$\\
 \hline
 $\delta_2(V_l)$&$0$&$2$&$(0,0)$&$(1,1)$\\
 \hline
   \end{tabular}
 \end{center}
 \end{Corollary}

\subsection{Reduction types ${\rm IV}^*,\;{\rm III}^*$ and ${\rm II}^*$}
\label{subsec:red IV*,III*,II*}

For reduction types ${\rm IV}^*$, ${\rm III}^*$ and ${\rm II}^*$, we will not use Tate's algorithm to write down explicit equations for the components in $\Phi_K^2(E)$ because the defining equations of the components in $\Phi_K^m(E),\;m\ge 2,$ are complicated. Instead, we make use of our knowledge of the reduction types of $E$ over totally ramified field extensions of $K$, and use the projection formula to exploit the symmetry of the graphs associated to $E^{min}_k$.

\begin{Lemma}
\label{lem:type over extensions}
Let $E/K$ be an elliptic curve. Let $K_m=K(t^{1/m}).$ If $E/K$ has reduction type ${\rm IV}^*$, then $E/K_3$ has good reduction if $\Char k\ne3$. If $E/K$ has reduction ${\rm III^*}$, then $E/K_2$ has reduction ${\rm I_0}^*$ and $E/K_3$ has reduction ${\rm III}$.
 \end{Lemma}
\begin{Proof}
Let $t_m=t^{1/m}$.
We recall from (\cite{Sil2}, Chapter IV, \S9) that $E/K$ has a minimal Weierstrass equation
\[y^2+a_1xy+a_3y=x^3+a_2x^2+a_4x+a_6\]
with $\nu(a_1)\ge1,\nu(a_2)\ge 2,\nu(a_3)=2,\nu(a_4)\ge3,\nu(a_6)\ge5$ if the reduction is ${\rm IV}^*$, and $\nu(a_1)\ge1,\nu(a_2)\ge 2,\nu(a_3)\ge3,\nu(a_4)=3,\nu(a_6)\ge5$ if the reduction is ${\rm III}^*$. Making substitutions $x\mapsto t_m^{2(m-1)}x$ and $y\mapsto t_m^{3(m-1)}y$, we will have a minimal Weierstrass equation for $E/K_m$ of the form
\[y^2+a_1 t_m^{-(m-1)}xy+t_m^{-3(m-1)}a_3y=x^3+a_2t_m^{-2(m-1)}x^2+t_m^{-4(m-1)}a_4x+t_m^{-6(m-1)}a_6\]
 with reduction type as specified in the statement of the lemma.
\end{Proof}

If $E$ has one of the reduction types ${\rm IV}^*$, ${\rm III}^*$ or ${\rm II}^*$, then $E^{min}_k$ is as in the figures below. Recall that the multiplicity-$1$ component $\Gamma_i$ corresponds to $i\in\Z/n\Z$ under the isomorphism $\delta_1$. The numbers on the irreducible components refer to the multiplicities. \\

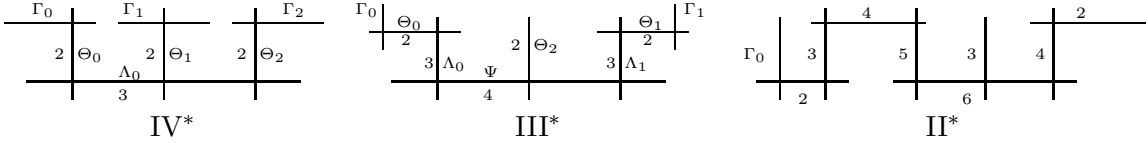
\setlength{\unitlength}{0.6cm}\linethickness{.2mm}
\begin{picture}(8,3)
{\tiny
\put(0,.7){\line(1,0){6}}\put(2,0.3){3}\put(2,0.8){$\Lambda_0$}
\put(1,0.3){\line(0,1){2}}\put(0.6,1.2){2}\put(1.1,1.2){$\Theta_0$}
\put(-0.5,2){\line(1,0){2}}\put(0.1,2.2){$\Gamma_0$}
\put(3,0.3){\line(0,1){2}}\put(2.6,1.2){2}\put(3.1,1.2){$\Theta_1$}
\put(2,2){\line(1,0){2}}\put(2.1,2.2){$\Gamma_1$}
\put(5,0.3){\line(0,1){2}}\put(4.6,1.2){2}\put(5.1,1.2){$\Theta_2$}
\put(4.5,2){\line(1,0){2}}\put(5.6,2.2){$\Gamma_2$}}
\put(2.7,-0.5){${\rm IV}^*$}
{\tiny
\put(8,0.7){\line(1,0){6}}\put(10,0.3){4}\put(10,0.8){$\Psi$}
\put(9,0.3){\line(0,1){1.8}}\put(9.1,1){$\Lambda_0$}\put(8.7,1){$3$}
\put(7.5,1.8){\line(1,0){2}}\put(8.1,1.9){$\Theta_0$}\put(8.2,1.5){2}
\put(7.8,1.4){\line(0,1){1}}\put(7.2,2.2){$\Gamma_0$}
\put(11,0.3){\line(0,1){1.8}}\put(10.6,1.4){2}\put(11.1,1.4){$\Theta_2$}
\put(13,0.3){\line(0,1){1.8}}\put(12.7,1){$3$}\put(13.1,1){$\Lambda_1$}
\put(12.5,1.8){\line(1,0){2}}\put(13.4,1.9){$\Theta_1$}\put(13.5,1.5){2}
\put(14.2,1.4){\line(0,1){1}}\put(14.4,2.2){$\Gamma_1$}}
\put(10.7,-0.5){${\rm III}^*$}
{\tiny
\put(16,.7){\line(1,0){2}}\put(16.9,0.2){$2$}
\put(16.5,0.3){\line(0,1){1.8}}\put(15.7,1.2){$\Gamma_0$}
\put(17.5,0.3){\line(0,1){2}}\put(17.1,1.2){$3$}
\put(17.2,2){\line(1,0){2.5}}\put(18.3,2.1){$4$}
\put(19,.7){\line(1,0){4}}\put(20.5,0.2){$6$}
\put(19.5,0.3){\line(0,1){2}}\put(19.1,1.2){$5$}
\put(21,0.3){\line(0,1){1.8}}\put(20.6,1.2){$3$}
\put(22.5,0.3){\line(0,1){2}}\put(22.1,1.2){$4$}
\put(22,2){\line(1,0){2.5}}\put(23,2.1){$2$}}
\put(19.7,-0.5){${\rm II}^*$}
\end{picture}\\

\begin{Proposition}
 \label{prop:sum of points IV*, III*, II*}
 Let $E/K$ be an elliptic curve with minimal proper regular model $E^{min}$. Let $\Xi_m$ denote a component of multiplicity-$m$ in $E^{min}_k$.
 \begin{myitemize}
   \item[(i)] If $E/K$ has reduction ${\rm IV}^*$ and $\Char k\ne 3$, then\\
\begin{center}
 \begin{tabular}{|c||c|c|c|c|}
  \hline
 $\Xi_m$& $\Theta_0$ & $\Theta_1$ & $\Theta_2$&$\Lambda_0$\\
 \hline
 $\delta_m(\Xi_m)$ & $0$ & $2$ & $1$ & $0$\\
  \hline
  \end{tabular}
 \end{center}
 \item[(ii)] If $E/K$ has reduction ${\rm III}^*$ and $\Char k\ne2,3$, then \\
 \begin{center}
 \begin{tabular}{|c||c|c|c|c|c|c|}
  \hline
 $\Xi_m$& $\Theta_0$ & $\Theta_1$ & $\Theta_2$&$\Lambda_0$&$\Lambda_1$&$\Psi$\\
 \hline
 $\delta_m(\Xi_m)$ & $0$ & $0$ & $1$ & $0$& $1$ & $0$\\
  \hline
  \end{tabular}
 \end{center}
 \item[(iii)] If $E/K$ has reduction ${\rm II}^*$, then $\delta_m(\Xi_m)=0$ for every irreducible component $\Xi_m$ in $E^{min}_k$.
  \end{myitemize}
 \end{Proposition}
 \begin{Proof}
 Let $K_m=K(t^{1/m}),\;m\in\{2,3,4\}.$ Let $\Gal(K_m/K)=\langle\sigma_m\rangle.$ We want to find $l\in\Phi_K(E)$ such that if $P\in E(K_m)$ reduces to a point on $\Xi_m$, then $\sum_{i=0}^{m-1}P^{\sigma_m^i}$ reduces to a point on $\Gamma_l$ and therefore $\delta_m(\Xi_m)=l.$ If $Q\in E(K)$, then we will denote the translation-by-$Q$ automorphism by $\tau_Q:E\to E.$ This $K$-automorphism extends to give an $\OK$-automorphism $\tau'_Q:E^{min}\to E^{min},$ see e.g. (\cite{Sil2}, Chapter IV, Proposition 4.6).
 In $(i)$ and $(ii)$, we know that $\Theta_0: y_1^2=0$ and $\delta_2(\Theta_0)=0$, see Lemma \ref{lem:y2=0, x3=0}.

$(i)$ Let $P\in E(K_2)$ be a point lying above $\Theta_i,\;i=1,2,$ and above no other component.
 Now let $Q\in E(K)$ be a point lying above $\Gamma_{2i}.$ We have $\tau'^*_Q(\Gamma_{i})=\Gamma_0,$ and the projection formula implies that $\tau'^*_Q(\Gamma_i).\tau'^*_Q(\Theta_i)=\Gamma_i.\Theta_i,$ see (\cite{Liubook}, Theorem 9.2.12), therefore $\tau'^*_Q(\Theta_i)=\Theta_0.$ Hence $\tau_Q(P)+\tau_Q(P)^{\sigma_2}\in E^0(K),$ because $\delta_2(\Theta_0)=0$, i.e., $(P+Q)+(P^{\sigma_2}+Q)\in E^0(K).$ In other words, $P+P^{\sigma_2}\in Q+E^0(K),$ i.e., $P+P^{\sigma_2}$ lies above $\Gamma_{2i}.$

 Since $E/K_3$ has reduction type ${\rm I_0}$, see Lemma \ref{lem:type over extensions}, $\kbar=k$ and $\Char k\ne3$, then $E(K_3)$ is divisible by $3$. For $P\in E(K_3)$ lying above $\Lambda_0$, there is $Q\in E(K_3)$ with $P=3Q$. Now $P+P^{\sigma_3}+P^{\sigma_3^2}=3(Q+Q^{\sigma_3}+Q^{\sigma_3^2})\in E^0(K)$ because $\Phi_E(K)\isom\Z/3\Z$.

 $(ii)$
 Following the same argument as in $(i)$ with $\tau_Q$ and $Q\in E(K)\setminus E^0(K)$, and remembering that $\delta_2(\Theta_0)=0$, we have
 $\delta_2(\Theta_1)=0$.
Let $\EE_m$ be the minimal proper regular model of $E/K_m$, $m=2,3$. Indeed, there is a morphism
$$ \psi_m:\EE_m\to\Norm(E^{min}\times_{\OK}\mathcal{O}_{K_m})\to E^{min},$$
where $\Norm(E^{min}\times_{\OK}\mathcal{O}_{K_m})$ is the normalisation of $E^{min}\times_{\OK}\mathcal{O}_{K_m}$, see for example (\cite{Liubook}, \S 10.4) and (\cite{Lorenzini}, pp. 10-11). The minimal proper regular model $\EE_m$ is obtained from the latter normalisation by contracting $(-1)$-curves.

By virtue of Lemma \ref{lem:type over extensions}, the reduction type of $E/K_2$ is ${\rm I}_0^*.$ There exists a minimal Weierstrass equation for $E/K_2$
 $$y^2+a_1xy+a_3y=x^3+a_2x^2+a_4x+a_6,\;a_i\in\mathcal{O}_{K_2},$$
 where $x^3+a_2x^2+a_4x+a_6=(x-t\alpha_1)(x-t\alpha_2)(x-t\alpha_3),\;\alpha_i\in\mathcal{O}_{K_2},$ such that $(\EE_2)_k$ consists of $\Gamma_0:z=0,$ $V_{-1}:y_1^2=0$ and the three distinct lines of $x_1^3+\tilde{a}_{2,1}x_1^2+\tilde{a}_{4,2}x_1+\tilde{a}_{6,3}=0.$ Let $P_i=(t\alpha_i,0).$ The reduction map $E(K_2)[2]=\{0,P_1,P_2,P_3\}\to \Phi_{K_2}(E)$ is surjective as each of these points lies on a different multiplicity-1 component. Now $\sigma_2$ fixes one of the $\alpha_i$'s and swaps the other two.
It follows that $\sigma_2$ fixes two multiplicity-1 components $\Gamma_{(0,0)}$ and $\Gamma'$ of $(\EE_2)_k$ and swaps the other two multiplicity-1 components $\Gamma_*,\Gamma^{\sigma_2}_*$. Note that $\Gamma_{(0,0)}$ and $\Gamma'$ lie above $\Gamma_0$ and $\Gamma_1$ in $E^{min}_k$ respectively.
 Since $\Char k\ne 2$ and $[K_2:K]$ divides the multiplicities of $\Theta_2$ and $\Psi,$ it is known that $\Gamma_*$ and $\Gamma^{\sigma_2}_*$ lie above $\Theta_2$ under the morphism $\psi_m$, see for example (\cite{Liubook}, Remark 10.4.8).
Therefore, if $P$ lies above $\Theta_2$ in $E^{min},$ then it lies above either $\Gamma_*$ or $\Gamma^{\sigma_2}_*$ in $\EE_2.$ Since $\Phi_{K_2}(E)\isom \mathbb{Z}/2\mathbb{Z}\times\mathbb{Z}/2\mathbb{Z},$ we have $P+P^{\sigma_2}$ lies above $\Gamma'$ in $(\EE_2)_k,$ and hence it lies above $\Gamma_1$ in $E^{min}_k.$
Furthermore, the multiplicity-$2$ irreducible component $V_{-1}$ of $(\EE_2)_k$ lies above $\Psi$. We notice that $K_2=K_4^G$ where $G=\{1,\sigma_4^2\}$. Therefore, if $P\in E(K_4)$ lies above $\Psi$, then $\sum_{i=0}^3P^{\sigma_4^i}=(P+P^{\sigma_4})+(P+P^{\sigma_4})^{\sigma_4^2}\in  E^0(K_2)$, because $\delta_2(V_{-1})=(0,0)$, i.e., this sum lies above $\Gamma_{(0,0)}$ of $(\EE_2)_k$ which lies above the identity component $\Gamma_0$ of $E^{min}_k$.

 We know that $E/K_3$ has reduction type ${\rm III}$, see Lemma \ref{lem:type over extensions}. If $P\in E(K_3)$, then $P$ lies on a multiplicity-$1$ component of $(\EE_3)_k$. Moreover, since $\Char k\ne 3$, according to the description of the morphism $\psi_3$ given in (\cite{Lorenzini}, pp. 10-11) if $P$ lies on $\Lambda_0$, then $P$ lies on the identity component of $(\EE_3)_k$, and if $P$ lies on $\Lambda_1$, then $P$ lies on the other multiplicity-$1$ component of $(\EE_3)_k$. Thus we are done by observing that $\Phi_{K_3}(E)\isom\Z/2\Z$.

  $(iii)$ This is straightforward because $\sum_{i=0}^{m-1}P^{\sigma^i_m}\in E(K),$ and $E(K)=E^0(K)$ for reduction type ${\rm II}^*.$
\end{Proof}

\section{Counting minimal models}
\label{sec:counting minimal degree n models}

Let $K$ be a Henselian discrete valuation field with ring of integers $\OK.$ We fix a uniformiser $t.$ We denote the normalised valuation on $K$ by $\nu.$ We assume that the residue field $k$ is perfect with algebraic closure $\kbar$.

Let $C$ be a smooth genus one curve over $K$ such that $C(K)\ne\emptyset.$
Let $E$ be the Jacobian elliptic curve of $C$ with minimal proper regular model $E^{min}$. When $E$ has additive reduction, we will assume moreover that $\Char k\ne 2$ for $n=2$, and $\Char k\ne2,3$ for $n=3,4.$

 Let $Q\in C(K)$. Again we fix an isomorphism $\beta:C\to E$ such that $\beta(Q)=0_E.$ So we can dispense with $C$ and write $E$ instead. If $D$ is an effective $K$-rational divisor of degree $n$ on $E,$ then the Riemann-Roch theorem implies that there exists a point $P\in E(K)$ such that $D\sim(n-1).0_E+P.$ Therefore, the equation defining $C\to\PP^{n-1}_K$ is an equation for the double cover of the projective line $E\to \PP_K^1$ given by the divisor class $[0_E+P]$ when $n=2,$ or the image of $E$ when it is embedded in $\mathbb{P}_K^{n-1}$ by the divisor class $[(n-1).0_E+P]$ when $n=3,4.$ Let  $\psi_E$ be the surjective homomorphism in the following short exact sequence
$$0\to E^0(K)\to E(K)\xrightarrow{\psi_E} \Phi_E(K)\to 0.$$

Let $S^n(P),\;n\le4$, be the set consisting of unordered tuples $(\Gamma_1,\ldots,\Gamma_l)$, where $\Gamma_i$ is an irreducible component of $E^{min}_{\kbar}$, satisfying the following conditions:
\begin{myitemize}
\item[(i)] $\sum_{i=1}^l\mult_{\kbar}(\Gamma_i)=n.$
\item[(ii)] $\sum_{i=1}^l\delta_{m_i}(\Gamma_i)=\psi_E(P)$, where $m_i=\mult_{\kbar}(\Gamma_i)$.
\item[(iii)] $(\Gamma_1^{\sigma},\ldots,\Gamma_l^{\sigma})=(\Gamma_1,\ldots,\Gamma_l)$ for every $\sigma\in\GG.$
\end{myitemize}
 We define the map $$\lambda_n:\{\textrm{minimal degree-$n$-models for $E\to\PP^{n-1}_K$ up to isomorphism}\}\to S^n(P);\;(\CC,\alpha)\mapsto (\Gamma_1,\ldots,\Gamma_l),$$ such that
\begin{myitemize}
\item[(i)] $(\Gamma_1,\ldots,\Gamma_l)$ consists of the irreducible components of the strict transform of $\CC_{\kbar}$ in $E^{min}_{\kbar}.$
\item[(ii)] If $\Gamma$ is the strict transform of a component $\Gamma'$ in $\CC_{\kbar}$, then $\Gamma$ appears in $(\Gamma_1,\ldots,\Gamma_l)$ as many times as $\deg_{\kbar}\Gamma'.$ In particular, we have $\sum_{i=1}^l\mult_{\kbar}(\Gamma_i)=n$. See \cite{Poonen1} and \cite{Bromwich} for the complete classification of the multiplicities and degrees of the irreducible components of $\CC_k$ when $n=3,4$. The case $n=2$ is clear.
\end{myitemize}

\begin{Lemma}
\label{lem:lambda is well defined}
The map $\lambda_n$ is well-defined.
\end{Lemma}
\begin{Proof}
To show that $\lambda_n$ is well-defined we need to prove (i) if $\CC$ and $\CC'$ are minimal isomorphic degree-$n$-models for $E\to\PP^{n-1}_K$, then $\lambda_n(\CC)=\lambda_n(\CC')$, and (ii) $\lambda_n(\CC)\in S^n(P)$. To prove (i) we know that the special fibers of any two isomorphic degree-$n$-models for $E\to\PP^{n-1}_K$ have the same irreducible components with the same types, therefore the corresponding tuples of both models are the same, see Theorem \ref{thm:injectivity2}.

 For (ii) we will show first that if $(\Gamma_1,\ldots,\Gamma_l)$ is the tuple associated to $\CC$, then $\sum_{i=1}^l\delta_{m_i}(\Gamma_i)=\psi_E(P)$ where $m_i=\mult_k(\Gamma_i).$ Let $D$ be a divisor on $E^{min}$ such that $\CC$ is obtained from $E^{min}$ by contraction using $D$, see Theorem \ref{thm:minimality}. The divisor $D$ intersects the $\Gamma_i$'s and no other components. Endowing the $\Gamma_j$'s with the reduced structure, we have $D. \Gamma_j=\deg_k\Gamma_j'$, where $\Gamma_j'$ is an irreducible component of $\CC_k$ whose strict transform in $E^{min}_k$ is $\Gamma_j$. Let $(x_1,\ldots,x_l)$ be a tuple of all intersection points of $D$ with the irreducible components $(\Gamma_1,\ldots,\Gamma_l)$ where $x_i\in \Gamma_i$. Note that $x_i$ may be repeated in $(x_1,\ldots,x_l)$ if $\Gamma_i$ is the strict transform of a component whose degree is greater than $1$. We have $D|_E\sim (n-1).0+P.$

Let $m_r=\max\{m_j:1\le j\le l\}$. It is clear that for $\Gamma_j,\;j=1,\ldots,l,$ the multiplicity $m_j\in\{1,m_r\}$, since otherwise we will have $1<m_j<m_r$ which implies that $m_r\ge 3$, hence $m_r+m_j\ge 5$ which contradicts that $\sum_{i=1}^l m_i=n\le 4$. By virtue of Proposition \ref{prop:auxiliary divisors}, there exists a closed point $P_j\in E$ such that $[K(P_j):K]=m_j$ and $\overline{\{P_j\}}\cap\Gamma_j=\{x_j\}$.  Let $K_j=K(P_j)$ and $\Gal(K_{j}^{un}/K^{un})=\langle\sigma_{j}\rangle.$ Let $\EE$ be the minimal proper regular model of $E/K_{r}^{un}$. Since $x_j$ lies on one and only one component of $E^{min}_{\kbar}$, there are exactly $m_j$ points $\{y_1,\ldots,y_{m_j}\}$ of $\EE$ lying above $x_j$, and each of these points lies on a multiplicity-$1$ component of $\EE_{\kbar}$, see (\cite{Liubook}, Remark 10.4.8). Now we view $D$ as a divisor on $\EE$. For $Q\in E$, $\tilde{Q}$ will denote its reduction. We have
$$\tilde{P}=(\summ D|_E)^{\sim}=\summ D|_{\EE_{\kbar}}=\sum_{j=1}^l\sum_{i=1}^{m_j}y_i=\sum_{j=1}^l\sum_{i=1}^{m_j}\tilde{P}_j^{\sigma_{j}^i}.$$
The second equality holds because $D$ intersects $\EE_{\kbar}$ only in multiplicity-$1$ components.
 Therefore, we have $P-\sum_{j=1}^l\sum_{i=1}^{m_j}P_j^{\sigma_{j}^i}\in E^0(K)$. Applying the surjective group homomorphism $\psi_E:E(K)\to\Phi_K(E)$, we get $\psi_E(P)=\sum_{j=1}^l\psi_E(\sum_{i=1}^{m_j}P_j^{\sigma_{j}^i})=
 \sum_{j=1}^l\delta_{m_j}(\Gamma_j).$

 To show that $\lambda_n(\CC)$ is $\GG$-invariant, we know that $D$ intersects the irreducible components of the strict transform of $\CC_{\kbar}$ in $E^{min}_{\kbar}$ and no other components. If $\Gamma$ is the strict transform of an irreducible component in $\CC_{\kbar}$, then there is an $x\in(D|_{E^{min}_{\kbar}})\cap\Gamma$. Moreover, $H^0(E^{min},\mathcal{O}_{E^{min}}(mD))$ is a free $\OK$-module for every $m\ge0$. Therefore, $D|_{E^{min}_{\Kbar}}$ and $D|_{E^{min}_{\kbar}}$ are $K$- and $k$-rational divisors respectively. Thus $x^{\sigma}\in(D|_{E^{min}_{\kbar}})^{\sigma}\cap\Gamma^{\sigma}=(D|_{E^{min}_{\kbar}})\cap\Gamma^{\sigma}$ for any $\sigma\in\GG$. In particular, $\Gamma^{\sigma}$ is the strict transform of an irreducible component of $\CC_{\kbar}$.
 \end{Proof}

We will show that $\lambda_n$ is bijective, and so the set of isomorphism classes of minimal degree-$n$-models for $E\to\PP^{n-1}_K$ is identified with $S^n(P)$. Consequently, we have the following theorem.

\newtheorem{thm2}{Theorem}[section]
\begin{Theorem}
\label{thm:counting models}
Let $E/K$ be an elliptic curve and $P \in E(K)$. Let
$E \to \PP^{n-1}_K$ be the morphism determined by the divisor class $[(n-1).0_E + P],\;n\in\{2,3,4\}$. Let $\delta_m : \Phi_K^m(E) \to \Phi_K(E)$ be the function defined in \S \ref{sec:computing in E min}.
Then there is a bijection between the set of minimal degree-$n$-models for $E\to \PP^{n-1}_K$ up to isomorphism
and the disjoint union of the following $\GG$-invariant sets.
\begin{myitemize}
\item[(i)] The set of unordered $n$-tuples
\[ S_1(n,P) = \{ (a_1,\ldots,a_n) : a_i \in \Phi_K(E) \mid a_1 + \ldots + a_n = \psi_E(P)\}^{\GG},\]
\item[(ii)] the set $S_n(P) = \{ a \in \Phi_K^n(E) \mid \delta_n(a) = \psi_E(P) \}^{\GG},$
\end{myitemize}
and

if $n\ge3$
\begin{myitemize}
\item[(iii)] the set $S_{(n-1,1)}(P) = \{ (a,b) \in \Phi_K^{n-1}(E) \times \Phi_K(E) \mid
\delta_{n-1}(a) + b = \psi_E(P) \}^{\GG},$
\end{myitemize}
and

if $n=4$
\begin{myitemize}
\item[(iv)] the set $S_{(2,1,1)}(P) = \{ (a,b,c) \in \Phi_K^{2}(E) \times \Phi_K(E)\times \Phi_K(E) \mid
\delta_{2}(a) + b + c = \psi_E(P) \}^{\GG},$

where the order of $b$ and $c$ is immaterial,
\item[(v)] the set of unordered pairs $$S_{(2,2)}(P) = \{ (a,b) \in \Phi_K^{2}(E) \times \Phi_K^2(E) \mid
\delta_{2}(a) + \delta_2(b) = \psi_E(P) \}^{\GG}.$$
\end{myitemize}
\end{Theorem}

  We will denote the number of multiplicity-1 irreducible components of $E^{min}_{\kbar}$ which are defined over $k$ by $c_k,$ this is the {\em Tamagawa number} of $E/K.$ Now we get the following corollary as a direct consequence of Theorem \ref{thm:counting models}.

\begin{Corollary}
\label{cor:counting reduced models}
Let $E/K$ be an elliptic curve and let $P \in E(K)$. Let
$E \to \PP^{n-1}_K$ be the morphism determined by the divisor class $[(n-1).0_E + P]$. The number of minimal degree-$n$-models for $E\to\PP^{n-1}_K$ is determined according to Table \ref{table1}. \\
 \begin{sidewaystable}
 \caption{}
 \centering
 \footnotesize{
 \begin{tabular}{|c|c||c|c|c|}
  \hline
 &$c_p$&$n=2$  & $n=3$ & $n=4$ \\
 \hline
 \hline
 \multirow{2}{*}{${\rm I}_{2m}$} &\multirow{2}{*}{$2$}& $m+1$, if $P\in E^0(K)$ &\multirow{2}{*}{$m+1$}& $(m+1)(m+2)/2$, if $P\in E^0(K)$\\
  && $1$, if $P\not\in E^0(K)$ & & $m+1$, if $P\not\in E^0(K)$ \\
 \hline
 \multirow{3}{*}{$\rm{I}_{2m}$}&\multirow{3}{*}{$2m$}&$m+1$ if $\psi_E(P)\in2\Phi_K(E)$&$(m+1)(2m+1)/3$ if $3\nmid m$&$m(m+1)(m+2)/3$ if $2\nmid \psi_E(P)$\\
 &&&$m(2m+3)/3+1$ if $3\mid\gcd(m,\psi_E(P))$&$m(m+1)(m+2)/3+1$ if $2\mid \psi_E(P)$\\
 &&$m$ if $\psi_E(P)\not\in2\Phi_K(E)$&$m(2m+3)/3$ if $3\mid m$ and $3\nmid \psi_E(P)$&$m(m+1)(m+2)/3+2$ if $4\mid \psi_E(P)$\\
 \hline
 ${\rm I}_{2m+1}$ &$1$& $m+1$ &$m+1$& $(m+1)(m+2)/2$ \\
 \hline
 \multirow{3}{*}{${\rm I}_{2m+1}$}&\multirow{3}{*}{$2m+1$}&\multirow{3}{*}{$m+1$}&$(m+1)(2m+3)/3$ if $3\nmid 2m+1$&\multirow{3}{*}{$(m+1)(m+2)(2m+3)/6$}\\
 &&&$(m+2)(2m+1)/3+1$ if $3\mid \gcd(2m+1,\psi_E(P))$&\\
 &&&$(m+2)(2m+1)/3$ if $3\mid 2m+1$ and $3\nmid\psi_E(P)$&\\
 \hline
 \hline
 ${\rm II}$&$1$&$1$&$1$&$1$\\
 \hline
 ${\rm III}$&$1$&$1$ if $\psi_E(P)=1$, $2$ if $\psi_E(P)=0$&$2$&$2$ if $\psi_E(P)=1$, $3$ if $\psi_E(P)=0$\\
 \hline
 \multirow{2}{*}{${\rm IV}$}&$1$&\multirow{2}{*}{$2$}&$2$&$3$\\
 &$3$&&$4$ if $\psi_E(P)=0$, $3$ otherwise &$5$\\
 \hline
 \hline
 $\rm{I}_0^*$&$1$&$2$&$3$&$4$\\
 \hline
  \multirow{2}{*}{${\rm I}_{2m}^*$}&\multirow{2}{*}{$2$}&$m+3$ if $P\in E^0(K)$&\multirow{2}{*}{$2m+4$}&$(m+2)(m+4)$ if $P\in E^0(K)$\\
 &&$m+2$ if $P\not\in E^0(K)$&&$(m+2)(m+3)$ if $P\not\in E^0(K)$\\
 \hline
 \multirow{3}{*}{${\rm I}_{2m}^*$}&\multirow{3}{*}{$4$}&$m+5$ if $\psi_E(P)=(0,0)$&\multirow{3}{*}{$2m+6$}&$(m+4)^2$ if $\psi_E(P)=(0,0)$\\
 &&$m+2$ if $\psi_E(P)=(1,1)$&&$(m+2)(m+5)$ if $\psi_E(P)=(1,1)$\\
 &&$2$ otherwise&&$4m+10$ otherwise\\
 \hline
 \multirow{2}{*}{${\rm I}_{2m+1}^*$}&\multirow{2}{*}{$2$}&$m+4$ if $P\in E^0(K)$&\multirow{2}{*}{$2m+5$}&$(m+3)(m+4)$ if $P\in E^0(K)$\\
 &&$m+2$ if $P\not\in E^0(K)$&&$(m+2)(m+4)$ if $P\not\in E^0(K)$\\
 \hline
 \multirow{3}{*}{${\rm I}_{2m+1}^*$}&\multirow{3}{*}{$4$}&$m+4$ if $2\mid\psi_E(P)$&\multirow{3}{*}{$2m+7$}&$(m+3)(m+6)$ if $\psi_E(P)=(0,0)$\\
 &&&&$(m+4)^2$ if $\psi_E(P)=2$\\
 &&$2$ otherwise&&$4m+12$ otherwise\\
 \hline
 \hline
 \multirow{2}{*}{${\rm IV}^*$}&$1$&\multirow{2}{*}{$3$}&$4$&$8$\\
 &$3$&&$8$ if $\psi_E(P)=0$, $6$ otherwise & $14$\\
 \hline
 ${\rm III}^*$&$2$& $4$ if $\psi_E(P)=0$, $2$ otherwise & $6$ & $15$ if $\psi_E(P)=0$, $10$ otherwise\\
 \hline
 ${\rm II}^*$&$1$& $3$ & $5$ & $10$\\
 \hline
   \end{tabular}}
   \label{table1}
   \end{sidewaystable}
\end{Corollary}
\begin{Proof}
For full details of the proof see \S\S 7,8 of \cite{SadekThesis}.
\end{Proof}

\section{Proof of Theorem \ref{thm:counting models}}
\label{sec:counting models}

We need a few lemmas before we proceed with the proof of Theorem \ref{thm:counting models}.

\begin{Lemma}
 \label{lem:divisibility of E0}
  Let $E/K$ be an elliptic curve with additive reduction. Assume that $\gcd(d,\Char k)=1$. Then the group $E^0(K)$ is divisible by $d$.
 \end{Lemma}
 \begin{Proof}
   Recall that $E^1(K)=\{P\in E(K):\tilde{P}=\tilde{0}_E\}.$ The group $E^0(K)/E^1(K)$ is isomorphic to $k^+$ because  $E$ has additive reduction. In particular, $E^0(K)/E^1(K)$ is divisible by $d$ because $(d,\Char k)=1.$ But by the theory of formal groups, the group $E^1(K)$ is uniquely divisible by $d$. Therefore, $E^0(K)$ is divisible by $d.$
 \end{Proof}

 If $E/K$ has non-split reduction type, then there exists a smallest finite unramified extension $L/K$ of degree $d$ over which $E$ has split reduction. In fact, $d=2$ except possibly when $E$ has non-split reduction of type ${\rm I}_0^*$ then $d\in\{2,3\}$.
 We define the norm map $\Norm_{L/K}$ to be $\Norm_{L/K}: E(L)\to E(K);\;Q\mapsto\sum_{i=1}^d Q^{\sigma^i}$,
where $\Gal(L/K)=\langle\sigma\rangle$.
\begin{Lemma}
\label{lem:trace is surjective}
Let $L$ be the smallest unramified field over which $E$ has split reduction. Then $\Norm_{L/K}:E^0(L)\to E^0(K)$ is surjective.
\end{Lemma}
\begin{Proof}
We first treat the additive case. So assume that $E$ has non-split additive reduction. Let $Q\in E^0(K).$ Since $E^0(K)$ is divisible by $d=[L:K]$,
see Lemma \ref{lem:divisibility of E0}, there is $Q'\in E^0(K)$ with $dQ'=Q.$ Now we have $\Norm_{L/K}(Q')=dQ'=Q.$

Now assume $E$ has non-split multiplicative reduction. We always have $[L:K]=2$. The non-singular reduction of $E$ will be denoted by $\tilde{E}_{ns}.$ Let $\ell$ be the residue field of $L.$ We denote the image of $\sigma$ under the isomorphism $\Gal(L/K)\isom\Gal(\ell/k)$ by $\sigma$ again. Let $\Norm_{\ell/k}:\tilde{E}_{ns}(\ell)\to\tilde{E}_{ns}(k);\;Q\mapsto Q+Q^{\sigma}$. Consider the following diagram.\\
\begin{displaymath}
\xymatrix{
0\ar[r]&E^1(L)\ar[d]_{\Norm_{L/K}}\ar[r]&E^0(L)\ar[d]_{\Norm_{L/K}}\ar[r]&\tilde{E}_{ns}(\ell)\ar[d]_{\Norm_{\ell/k}}\ar[r]&0\\
0\ar[r]&E^1(K)\ar[r]&E^0(K)\ar[r]&\tilde{E}_{ns}(k)\ar[r]&0
 }
 \end{displaymath}

To prove that $\Norm_{L/K}:E^0(L)\to E^0(K)$ is surjective, we only need to show the surjectivity of both $\Norm_{L/K}:E^1(L)\to E^1(K)$ and $\Norm_{\ell/k}:\tilde{E}_{ns}(\ell)\to\tilde{E}_{ns}(k)$. Let $Q\in E^1(K).$ Since $E^1(K)$ is divisible by $2,$ there is a $Q'\in E^1(K)$ such that $2Q'=Q.$ So, $\Norm_{L/K}(Q')=2Q'=Q.$

Now we will show the surjectivity of $\Norm_{\ell/k}.$ The isomorphism $f:\tilde{E}_{ns}(\ell)\isom \ell^*$ induces an isomorphism $\tilde{E}_{ns}(k)\isom U:=\{u\in \ell^*:\Norm_{\ell/k}(u)=1\}$. Moreover, the $\ell$-automorphism $\sigma^{-1}f\sigma f^{-1}:\ell^*\to \ell^*$ is $u\mapsto u^{-1},$ see (\cite{Liubook}, Exercise 10.2.7).
Therefore, the map $\Norm_{\ell/k}$ induces the map $\ell^*\to U;\;v\mapsto v/v^{\sigma}.$
Now according to Hilbert's Theorem 90, for $u\in \ell^*$ we have $\Norm_{\ell/k}(u)=1$ if and only if $u=v/v^{\sigma}$
for some $v\in \ell^*.$ Therefore, $\Norm_{\ell/k}$ is surjective.
\end{Proof}

\begin{Lemma}
 \label{lem:constructing divisors}
 Let $E/K$ be an elliptic curve with minimal proper regular model $E^{min}.$ Let $P\in E(K)$ and $(\Gamma_1,\ldots,\Gamma_l)$ be an element of $S^n(P)$ endowed with the reduced structure.
 Then there exists a divisor $D$ on $E^{min}\to\Spec\mathcal{O}_{K^{un}}$ such that:
  \begin{myitemize}
  \item[(i)] $(D|_E)^{\sigma}=D|_E$ for every $\sigma\in\Gal(K^{un}/K).$
  \item[(ii)] $D. \Gamma_i=d_i,$ where $d_i$ is the number of times $\Gamma_i$ appears in $(\Gamma_1,\ldots,\Gamma_l)$.
  \item[(iii)] $D|_E\sim(n-1).0+P.$
  \end{myitemize}
 \end{Lemma}
 \begin{Proof}
 Assume first that $E$ has reduction of type ${\rm I}_m,m\ge0$. Then $E^{min}_{\kbar}$ contains only multiplicity-$1$ irreducible components and we have $l=n$. Let $\ell=k(\Gamma_1,\ldots,\Gamma_n)$. Let $L/K$ be the unramified extension with residue field $\ell/k$ and $\Gal(L/K)=\langle\sigma\rangle$. We have $[L:K]=[\ell:k]\le 2.$ Now we pick a point above each $\Gamma_i$. If $\Gamma_i$ is defined over $k$, then choose $Q_i\in E(K)$ above $\Gamma_i$. If $\Gamma_i$ is defined over $\ell$, then pick $Q_i\in E(L)$ lying above $\Gamma_i$. Since $(\Gamma_1,\ldots,\Gamma_n)$ is $\GG$-invariant, then $\Gamma_i^{\sigma}=\Gamma_j$ for some $j\in\{1,\ldots,n\}$. We pick $Q_j:=Q_i^{\sigma}$ to be the point above $\Gamma_j$. Since $\sum_{i=1}^n\delta_{1}(\Gamma_i)=\psi_E(P)$, it follows that $P':=\sum_{i=1}^nQ_i-P\in E^0(K).$ If there exists $i$ such that $Q_i\in E(K)$, then replace $Q_i$ with $Q_i-P'$. Otherwise, there exists $i,j$ such that $Q_i\in E(L)\setminus E(K)$ and $Q_j=Q_i^{\sigma}$.
According to Lemma \ref{lem:trace is surjective}, there is $Q'\in E^0(L)$ with $\Norm_{L/K}(Q')=P'$. Now replace $Q_i$ with $Q_i-Q'$, and $Q_j$ with $Q_j-Q'^{\sigma}$. We set our divisor to be $D=\sum_{i=1}^n\overline{\{Q_i\}}$.

 Now assume that $E$ has additive reduction. If $\Gamma_i$ is a multiplicity-$1$ component, then we pick $Q_i$ above $\Gamma_i$ as in the multiplicative case. We note that for non-split ${\rm I}_0^*$, $\Gamma_i$ might be defined over a cubic extension of $k$. In case $\Gamma_i$ is defined over $k$ and of multiplicity $m_i\ge2$, then we choose $Q_i\in E$ above $\Gamma_i$, where $K(Q_i)/K$ is a totally ramified extension with $[K(Q_i):K]=m_i$. In case of non-split reduction of type ${\rm IV}^*$ and $n=4$, we can have $l=2$ where $\Gamma_1$ is of multiplicity-$2$ and defined over a quadratic extension $\ell/k$. Then we pick $Q_1$ above $\Gamma_1$ where $[L(Q_1):L]=2$ and we choose $Q_2:=Q_1^{\sigma}$.
Let $G(Q_i)=\sum_{j=1}^{m_i}Q_i^{\sigma_i^j}$ where $\Gal(K(Q_i)^{un}/K^{un})=\langle\sigma_i\rangle$. Since $\sum_{i=1}^l\delta_{m_i}(\Gamma_i)=\psi_E(P)$, it follows from the definition of $\delta_{m_i}$ that $P':=\sum_{i}G(Q_i)-P\in E^0(K).$
By virtue of Lemma \ref{lem:divisibility of E0}, there is $Q'\in E^0(K)$ with $nQ'=P'$. We set $Q'_i=Q_i-Q'$ for every $i$, and $D=\sum_{i}\overline{\{Q_i'\}}$.
 \end{Proof}

\begin{ProofOf}{Theorem \ref{thm:counting models}} We need to prove that the well-defined map $$\lambda_n:\{\textrm{ minimal degree-}n\textrm{-models for }E\to\PP^{n-1}_K\textrm{ up to isomorphism}\}\longrightarrow S^n(P)$$ is bijective.
To prove that $\lambda_n$ is injective, let $(\CC_1,\alpha_1)$ and $(\CC_2,\alpha_2)$ be two minimal degree-$n$-models for $E\to\PP^{n-1}_K$. Assume that $\lambda_n(\CC_1)=\lambda_n(\CC_2).$ We want to show that $(\CC_1,\alpha_1)$ and $(\CC_2,\alpha_2)$ are isomorphic. The fact that they have the same corresponding tuples implies that they have the same strict transform in $E^{min}_k$, therefore $(\CC_1,\alpha_1)$ and $(\CC_2,\alpha_2)$ are isomorphic, see Theorem \ref{thm:injectivity2}.

Now we will prove that $\lambda_n$ is surjective. So assume that $(\Gamma_1,\ldots,\Gamma_l)\in S^n(P)$ and we want to construct a minimal degree-$n$-model $(\CC,\alpha)$ whose image under $\lambda_n$ is $(\Gamma_1,\ldots,\Gamma_l)$. By virtue of Lemma \ref{lem:constructing divisors}, there exists a divisor $D$ on $E^{min}$ such that $D|_E$ is $K$-rational, $D$ intersects $\Gamma_i$ as many times as its occurrence in the tuple, and $D|_E\sim(n-1).0+P.$  Consider the following $S$-model for $E$
$$\CC:= \Proj(\bigoplus_{m=0}^{\infty}H^0(E^{min},\mathcal{O}_{E^{min}}(mD))).$$
The model $\CC$ is a minimal degree-$n$-model for $E\to\PP^{n-1}_K$, see \cite{SadekMinimal}. The special fiber $\CC_k$ consists of the components in $(\Gamma_1,\ldots,\Gamma_l)$ because $D$ intersects these components in $E^{min}_k$ and intersects no other components. Moreover, each component of $\CC_k$ has degree equal to the number of iterations of its strict transform in $(\Gamma_1,\ldots,\Gamma_l)$. We obtain $\alpha$ from the linear equivalence of $D|_E$ and $(n-1).0_E + P.$
\end{ProofOf}

 \section{Counting minimal global models}
\label{sec:counting global models}

In this section we will be interested in attacking the global question. Let $K$ be a number field of class number 1 with ring of integers $\OK$. If $\pp$ is a prime in $K$, then we will denote the completion of $K$ at $\pp$ by $K_{\pp}$, and its ring of integers by $\mathcal{O}_{K_{\pp}}$. Let $C$ be a smooth genus one curve defined by an integral genus one equation of degree $n$ over $K$.  We assume moreover that $C(K_{\pp})\ne\emptyset$ for every prime $\pp$ of $K$. A {\em minimal global degree-$n$-model $(\CC,\alpha)$ for $C\to\PP^{n-1}_{K}$} consists of an $\Spec\OK$-scheme defined by an $\OK$-integral genus one equation of degree $n$ which is $\mathcal{O}_{K_{\pp}}$-minimal at every prime $\pp$, and an isomorphism $\alpha:\CC_K\isom C$ defined by an element in $\mathcal{G}_n(K)$. An {\em isomorphism} $\beta:(\CC_1,\alpha_1)\isom(\CC_2,\alpha_2)$ of degree-$n$-models for $C\to\PP^{n-1}_K$ is an isomorphism $\beta:\CC_1\isom\CC_2$ of $\Spec\OK$-schemes defined by an element in $\mathcal{G}_n(\OK)$ with $\beta_K=\alpha_2^{-1}\alpha_1$.

It is known that a minimal global degree-$1$-model for $C/K$, i.e., a Weierstrass model, exists and is unique up to isomorphism. A proof of the existence of a minimal global degree-$n$-model for $C/K$, where $C(K_{\pp})\ne\emptyset$ for every $\pp$, when $n\in\{2,3,4\}$ can be found in \cite{FiStCr}. The author gave a new proof in \cite{SadekMinimal}.

The problem of counting minimal global degree-$n$-models up to $\OK$-isomorphism can be tackled locally. We consider $C$ as a curve over $K_{\pp}$ and use the previous section to determine the number $N_{\pp}$ of minimal degree-$n$-models for $C\to\PP_{K_{\pp}}^{n-1}$ up to $ \mathcal{O}_{K_{\pp}}$-isomorphism. It turns out that we need only to investigate a finite set of primes. Then we collect the local data at the finite places using strong approximation.

If $m\in\OK,$ then we set $P(m)=\{\pp:\pp\textrm{ is prime, }\mathfrak{p}^2\mid m\}.$
The main result of this section is the following Theorem.

\begin{Theorem}
\label{thm:global models}
Let $K$ be a number field of class number $1$. Let $C\to\PP^{n-1}_{K}$ be a smooth curve defined by a genus one equation of degree $n,\;n=2,3,4.$ Assume that $C(K_{\pp})\ne\emptyset$ for every prime $\pp$. Let $E/K$ be the Jacobian elliptic curve of $C$ with minimal discriminant $\Delta$. Let $N$ and $N_{\pp}$ denote the numbers of minimal global degree-$n$-models for $C\to\PP^{n-1}_{K},$ up to $\OK$-isomorphism, and minimal degree-$n$-models for $C\to\PP^{n-1}_{K_{\pp}},$ up to $\mathcal{O}_{K_{\pp}}$-isomorphism, respectively. Then
$$N=\prod_{\pp\in P(\Delta)}N_{\pp}.$$
\end{Theorem}

Theorem \ref{thm:global models} follows immediately from the bijectivity of the following map
\begin{eqnarray}
\{\textrm{minimal global degree-}n\textrm{-models for } C/K\}&\longrightarrow&\prod_{\pp\in P(\Delta)}\{\textrm{minimal degree-}n\textrm{-models for }C/K_{\pp}\}\nonumber\\
(\CC,\alpha)&\mapsto&((\CC,\alpha),\ldots,(\CC,\alpha)).\nonumber
\end{eqnarray}
 Notice that the above two sets of degree-$n$-models are defined up to isomorphism.
Furthermore, our work always enables us to compute $N_{\pp}$ when $\pp$ lies above a prime $p\ge5.$
Before proceeding with the proof of Theorem \ref{thm:global models} we need the following Lemmas.
\begin{Lemma}
\label{lem:matrices}
Let $A\in\GL_n(K_{\pp})\cap\Mat_n(\mathcal{O}_{K_{\pp}})$ have coprime entries. Assume that $\pp=\pi\OK$. Then there exist matrices $U\in \GL_n(\OK)$ and $V\in\GL_n(\mathcal{O}_{K_{\pp}})$ such that $A=VDU,$ where $D=\diag(\pi^{r_1},\ldots,\pi^{r_{n-1}},1)$ and $r_1\ge\ldots\ge  r_{n-1}.$
\end{Lemma}
\begin{Proof}
We claim that there exists a matrix $B\in\GL_n(K)\cap\Mat_n(\OK)$ such that $V':=BA^{-1}\in\GL_n(\mathcal{O}_{K_{\pp}}).$ Granted this claim we write the Smith Normal Form for the matrix $B,$ so we have $B=U'D'DU$ where $U,U'\in\GL_n(\OK),$ $D'$ is a diagonal matrix whose entries are not divisible by $\pi$, and $D=\diag(\pi^{r_1},\ldots,\pi^{r_{n-1}},1),\;r_1\ge\ldots\ge  r_{n-1}.$ Then we set $V:=V'^{-1}U'D'\in\GL_n(\mathcal{O}_{K_{\pp}}),$ hence we are done.

To prove the claim, assume that $A=(a_{ij})_{i,j}.$ Recall that every element in $\mathcal{O}_{K_{\pp}}$ can be written uniquely as $\sum_{i\ge0} a_i \pi^i,$ where $a_i\in\OK$ lies in a set of representatives of $\mathcal{O}_{K}/\pp.$ By a continuity argument, there exists an integer $m>0$ such that the matrix $B=(a_{ij}\mod \pi^m)_{i,j}\in\GL_n(K)\cap\Mat_n(\OK)$. Now we have $BA^{-1} \equiv \id_n\mod \pi^m$ and hence $BA^{-1}\in\GL_n(\mathcal{O}_{K_{\pp}})$.
\end{Proof}

\begin{Lemma}
\label{lem:Zp isom Z}
Let $\phi$ and $\phi'$ be two minimal $\mathcal{G}_n(K_{\pp})$-equivalent genus one equations of degree $n$ with coefficients in $\OK$ and $\mathcal{O}_{K_{\pp}}$ respectively. Then $\phi'$ is $\mathcal{G}_n(\mathcal{O}_{K_{\pp}})$-equivalent to a minimal genus one equation of degree $n$ whose coefficients lie in $\OK$.
\end{Lemma}
\begin{Proof}
Assume that $\phi'$ is obtained from $\phi$ via $[\mu_n,A_n]$ in $\mathcal{G}_n(K_{\pp}).$ For $r\in K_{\pp}^*,$ the following transformations are identical:
$$[\mu_n,A_n]=[r^{-2}\mu_n,rA_n]\textrm{ when }n=2,4,\textrm{ and } [\mu_3,A_3]=[r^{-3}\mu_3,rA_3].$$
Therefore, we can assume that $A_n$ has coprime entries in $\mathcal{O}_{K_{\pp}}.$ Lemma \ref{lem:matrices} allows us to write $A_n=V_nD_nU_n$ where $V_n\in\GL_n(\mathcal{O}_{K_{\pp}})$, $U_n\in\GL_n(\OK)$, and $D_n=\diag(\pi^{r_1},\ldots,\pi^{r_{n-1}},1)$ where $\pp=\pi\OK$. Similarly, we can write $\mu_4=\nu_4'\tau_4\nu_4$ where $\nu_4'\in\GL_2(\mathcal{O}_{K_{\pp}})$, $\nu_4\in\GL_2(\OK)$ and $\tau_4=\diag(\pi^{-m_1},\pi^{-m_2}).$

Let $\psi$ be the $\OK$-integral genus one equation obtained from $\phi$ via the transformation $[1,U_n]$ when $n=2,3$, and via $[\nu_4,U_4]$ when $n=4$. Then $\psi$ lies in the same $\mathcal{G}_n(\OK)$-equivalence class as $\phi$. Let $\phi''$ be the genus one equation obtained from $\psi$ via the transformation $[\mu_n', D_n]$, where $\mu_2'=(\det D_2)^{-2}$,
$\mu_3'=(\det D_3)^{-1}$ and $\mu_4'=\tau_4$. It is clear that $\phi''$ is $\mathcal{G}_n(\mathcal{O}_{K_{\pp}})$-equivalent to $\phi'.$ We claim that $\phi''$ has coefficients in $\OK$. If it is not the case, then some of the coefficients of the polynomials defining $\phi''$ would lie in $\frac{1}{\pi^r}\OK$ for some $r>0$. But $\phi''$ is obtained from $\phi'$ via $[\omega_n,V_n^{-1}]$, where $\omega_n\in\mathcal{O}_{K_{\pp}}^*$ when $n=2,3$, and $\omega_4\in\GL_2(\mathcal{O}_{K_{\pp}})$. Since $\phi'$ is $\mathcal{O}_{K_{\pp}}$-integral, it follows that $\phi''$ should be $\mathcal{O}_{K_{\pp}}$-integral, which is a contradiction.
\end{Proof}

 The following lemma, (\cite{Sil2}, Chapter IV, Lemma 9.5), will be used to justify our choice of the set of prime numbers $P(\Delta).$

\begin{Lemma}
\label{lem:P(Delta)}
Let $E/K_{\pp}$ be an elliptic curve with discriminant $\Delta.$ If $\nu_{\pp}(\Delta)=1,$ then $E$ has reduction of type ${\rm I}_1.$
\end{Lemma}

\begin{ProofOf}{Theorem \ref{thm:global models}}
First we will show that the map $\lambda$ is well defined. Let $(\CC_1,\alpha_1)$ and $(\CC_2,\alpha_2)$ be two isomorphic minimal global degree-$n$-models for $C\to\PP^{n-1}_{K}.$ Then $\alpha:=\alpha_2^{-1}\alpha_1:(\CC_1)_{K_{\pp}}\to (\CC_2)_{K_{\pp}}$ is defined by an element in $\mathcal{G}_n(\OK)\hookrightarrow\mathcal{G}_n(\mathcal{O}_{K_{\pp}})$ for every prime $\pp,$ i.e., $(\CC_1,\alpha_1)$ and $(\CC_2,\alpha_2)$ have the same image under $\lambda.$

To show that $\lambda$ is injective, let $(\CC_1,\alpha_1)$ and $(\CC_2,\alpha_2)$ be two minimal global degree-$n$-models for $C\to\PP^{n-1}_{K}$ with the same image under $\lambda.$ We need to show that $(\CC_1,\alpha_1)$ and $(\CC_2,\alpha_2)$ are isomorphic. Let $\alpha:=\alpha_2^{-1}\alpha_1.$ The map $\alpha$ is defined by an element $[\mu,A]\in\mathcal{G}_n(K).$ We can assume that $A\in\Mat_n(\OK)$ has coprime entries. Since $(\CC_1,\alpha_1)$ and $(\CC_2,\alpha_2)$ have the same image under $\lambda,$ in particular they are $\mathcal{O}_{K_{\pp}}$-isomorphic for every $\pp\in P(\Delta)$, it follows that $A\in\GL_n(\mathcal{O}_{K_{\pp}}),$ and so $\pp\nmid \det A.$ If $\pp\not\in P(\Delta),$ then $E/K_{\pp}$ has either reduction types ${\rm I}_0$ or ${\rm I}_1,$ see Lemma \ref{lem:P(Delta)}. But according to Corollary \ref{cor:unique model when IO or I1}, when $E$ has either reduction types ${\rm I}_0$ or ${\rm I}_1,$ there is a unique degree-$n$-model for $C\to\PP^{n-1}_{K_{\pp}}.$ That means that for $\pp\not\in P(\Delta),$ $(\CC_1,\alpha_1)$ and $(\CC_2,\alpha_2)$ are isomorphic as degree-$n$-models for $C\to\PP^{n-1}_{K_{\pp}}.$ Hence $A\in \GL_n(\mathcal{O}_{K_{\pp}})$ for every prime $\pp,$ in particular $\pp\nmid\det A.$ Thus $\det A\in\mathcal{O}_K^*$ and $A\in\GL_n(\OK).$

Now we will prove the surjectivity of $\lambda.$ We will assume without loss of generality that the defining genus one equation $\phi$ of $C$ has coefficients in $\OK$ and that the associated discriminant is everywhere minimal.

If $P(\Delta)=\emptyset$, then the scheme defined by $\phi$ is the unique minimal global degree-$n$-model for $C\to\PP^{n-1}_K.$ So, assume $P(\Delta)=\{\pp_1,\ldots,\pp_m\},\;m\ge1$. Let $(\CC_i,\alpha_i),$ where $1\le i\le m,$ be a minimal degree-$n$-model for $C\to\PP^{n-1}_{K_{\pp_i}}.$ We want to construct a minimal global degree-$n$-model $(\CC,\alpha)$ for $C\to\PP^{n-1}_{K}$ such that $\alpha^{-1}\alpha_i:(\CC_i)_K\to\CC_K$ is defined by an element in $\mathcal{G}_n(\mathcal{O}_{K_{\pp_i}})$ for each $i.$ Let $\phi_i$ be the defining genus one equation of $\CC_i$.
By virtue of Lemma \ref{lem:Zp isom Z}, we can assume that $\phi_i$ has coefficients in $\OK$ and is obtained from $\phi$ via $[\mu_i,D_iU_i]$ where $D_i=\diag(\pi_i^{r_{i,1}},\ldots,\pi_i^{r_{i,n-1}},1)$, where $\pp_i=\pi_i\OK$, $U_i\in\GL_n(\OK)$ and $\mu_i$ is a scaling element. In fact, we can assume that $U_i\in\SL_n(\OK)$ by acting by a diagonal matrix whose entries lie in $\mathcal{O}_K^*$.
Given integers $m_i>0$, there exists a matrix $U\in\SL_n(\OK)$ such that $U\equiv U_i\mod \pi_i^{m_i}$ for every $i$, see (\cite{Fiminimise}, Lemma 3.2). We note that $\prod_{i=1}^m D_iUU_j^{-1}D_j^{-1}\equiv\prod_{i\ne j}D_i \mod \pi_j^{m_j}.$ Therefore, the genus one equation $\psi$ obtained from $\phi$ via the transformation $[\prod_{i=1}^m\mu_i,\prod_{i=1}^m D_iU]$ is $\mathcal{G}_n(\mathcal{O}_{K_{\pp_j}})$-equivalent to $\phi_j$ for every $j$.

 Now we want to show that $\psi$ is $\OK$-integral. This will imply that $\psi$ defines a minimal global degree-$n$-model for $C\to\PP^{n-1}_{K}$ which is $\mathcal{O}_{K_{\pp_j}}$-isomorphic to $(\CC_j,\alpha_j)$ for every $j$. Hence we will be done with the surjectivity. Assume on the contrary that $\psi$ is not $\OK$-integral. Since $\psi$ is obtained from the $\OK$-integral genus one equation $\phi$ via $[\prod_{i=1}^m\mu_i,\prod_{i=1}^m D_iU]$, it follows that some of the coefficients of the defining polynomials of $\psi$ lie in $\frac{1}{b}\OK,\;b=\pi_1^{l_1}\pi_2^{l_2}\ldots \pi_m^{l_m}$, where $l_i\ge 0$ and $l_j>0$ for some $j\in\{1,\ldots,m\}$. We have shown that $\psi$ is $\mathcal{G}_n(\mathcal{O}_{K_{\pp_j}})$-equivalent to the $\OK$-integral genus one equation $\phi_j$. It follows that $\psi$ is $\mathcal{O}_{K_{\pp_j}}$-integral, which is a contradiction.
\end{ProofOf}

The calculations included in the following examples are performed using \textsf{MAGMA}, see \cite{magma}. $N_n(T)$ will denote the number of minimal degree-$n$-models when the reduction of the Jacobian is of type $T$.

\begin{Example}
Consider the elliptic curve $E:y^2+xy=x^3-x^2-617x+5916/\Q$. It has reduction of types ${\rm III}^*$ and ${\rm I}_2$ at its bad primes $5$ and $19$ respectively. The following minimal global genus one equation $\phi_3$ of degree $3$ defines an everywhere locally soluble element $C$ in the $3$-Selmer group of $E$.
\begin{eqnarray}
21686353648850 x^3&+&234081254700017 x^2y+9338329782950x^2z+842219868972245xy^2+\nonumber\\
67198263238095 xyz&+&1340388284750 xz^2+1010096983050575 y^3+120889031707155 y^2z+\nonumber\\
4822691362750 yz^2&+&64131409475 z^3 =0. \nonumber
\end{eqnarray}
The number of minimal global degree-$3$-models for $C\to\PP^2_K$ is $N_3({\rm III}^*)\times N_3({\rm I}_2)=6\times 2=12$, see Theorem \ref{thm:global models}. These models have defining genus one equations obtained from $\phi_3$ via the following transformations in $\mathcal{G}_3(\Q)$.\\
\begin{tabular}{cccc}
$[1,\id_3]$, & $[1/5,\diag(5,1,1)]$, &  $[1/5,\diag(1,5,1)]$, &$[1/25,\diag(5,5,1)]$,\\
  $[1/25,\diag(5,1,5)]$, & $[1/25,\diag(1,25,1)]$, & $[1/19,\diag(1,1,19)]$, & $[1/95,\diag(5,1,19)]$,\\
 $[1/95,\diag(1,5,19)]$,& $[1/475,\diag(5,5,19)]$, & $[1/475,\diag(5,1,95)]$, & $[1/475,\diag(1,25,19)]$.
\end{tabular}\\

We notice that the genus one equations $\phi_3$ is not reduced in the sense of \cite{FiStCr}, but we preferred to write it in that way to make the transformations as simple as possible.
\end{Example}
\begin{Example}
Let $E:y^2+xy+y=x^3-4x-3/\Q$. The curve $C:y^2=-3x^4+2x^3+7x^2-2x-3$ represents an element in the $2$-Selmer group of $E$. A second $2$-descent on $C$ gives the following minimal global genus one equation $\phi_4$ of degree $4$.
\begin{eqnarray}
x_1^2-x_1x_3-x_2^2+x_2x_4+x_3^2&=&0,\nonumber\\
x_1x_4+x_2^2+x_2x_3-x_2x_4+x_3^2-x_3x_4&=&0.\nonumber
\end{eqnarray}
The equation $\phi_4$ defines a smooth genus one curve $C_4/\Q$.
The discriminant $\Delta$ of $E$ is $185$. Therefore, we have $P(\Delta)=\emptyset$.
Hence the minimal global degree-$4$-model $\CC$ defined by $\phi_4$ is unique, see Theorem \ref{thm:global models}.
\end{Example}

\hskip-18pt\emph{\bf{Acknowledgements.}}
This paper is based on the author's Ph.D. thesis \cite{SadekThesis} at Cambridge University. The author would like to thank his supervisor Dr. Tom Fisher for all the encouragement and guidance through this work.

\bibliographystyle{plain}
\footnotesize
\bibliography{thesisreferences}

\end{document}